# A CONSTRUCTIVE PROOF OF BEAL'S CONJECTURE

Nicholas J. Daras

*Hellenic Military Academy,*
*Department of Mathematics & Engineering Sciences,*
*16673 Vari Attikis, Greece*

E-mail: njdaras@sse.gr

**Abstract.** *We prove that there is no non-trivial integral positive solution to the generalized Fermat equation.*

In the early twentieth century, David Hilbert presented twenty-three great mathematical problems which featured main directions of scientific research throughout the period that followed ([3, 4]). By analogy, in 2016, John F. Nash and Michael Th. Rassias gave a list with seventeen, currently unsolved problems in modern mathematics, in the belief that these problems are expected to determine several of the main research directions at least during the beginning of the XXI century ([7]). The third problem in the series of this list refers to the exploration of the possibility of extending Fermat's Last Theorem.

Even before Andrew Wiles announced his famous proof of this Theorem ([8]), various generalizations had already been considered, to equations of the shape $Ax^p + By^q = Cz^r$, for fixed integers $A$, $B$ and $C$. In this direction, the Theorem of Henri Darmon and Andrew Granville states that *if $A, B, C, p, q$ and $r$ are fixed positive integers with $p^{-1} + q^{-1} + r^{-1} < 1$, the equation $Ax^p + By^q = Cz^r$ has at most finitely many solutions in coprime non-zero integers $x$; $y$ and $z$* ([2]). However, as is made clear through those mentioned by Michael A. Bennett, Imin Chen, Sander R. Dahmen and Soroosh Yazdani ([1, 7]), except the solutions identified by Preda Mihăilescu in the *Catalan equation* ([6]) and the solutions derived from some elementary numerical identities, in most cases, there is no non-trivial solution of this equation once we assume that A = B = C = 1.

In this connection, *Beal's Conjecture* argues that *if $\min\{p, q, r\} \geq 3$, there are no non-trivial co-prime integral solutions to the generalized Fermat equation $x^p + y^q = z^r$* ([2, 7]). So far, many computational attempts have produced strong indications that this conjecture may be correct ([5]). The aim of this paper is to give a proof of this Conjecture.

To this end, let us suppose that, on the contrary, Beal's Conjecture does not apply. This is equivalent to assuming that there are integer exponents $p \geq 3$, $q \geq 3$ and $r \geq 3$ and three **coprime** positive integers $x = even$, $y = odd$ and $z = odd$ satisfying a generalized Fermat equation of the form $x^p \pm y^q = z^r$.

Given any $m, n, \lambda$ in $\mathbb{Z}^+$, $\lambda = odd$, and $\sigma \in \mathbb{Z}^+ \setminus \{1\}$, put
$$\tau_{m,n} := -y^{\sigma q - 4} n(2m + n) + \tau_0, \text{ with } \tau_0 := \lambda y^{q-2}, \lambda = odd.$$

By adopting this definition, it is clear that $\tau_{m,n}$ is an integer and the following apply.

**Proposition 1**. *There exist $n^*, m^* \in \mathbb{Z}^+$, for which the corresponding expression $z^2 \mp y^2 \tau_{m^*,n^*}$ represents a prime number.*

*Proof.* Since $gcd(\mp y^{\sigma q-2}, \mp[z^2 \mp \tau_0 y^2]) = 1$, an application of Dirichlet's Theorem on arithmetic progressions, in its basic form, shows that there are infinitely many primes such that $\pi \equiv \mp(z^2 \mp \tau_0 y^2) \mod (\mp y^{\sigma q-2})$. In particular, there exist infinitely many natural numbers $X$ such that the integers $\mp y^{\sigma q-2} X + (z^2 \mp \tau_0 y^2)$ are primes $\pi$. Since $\tau_0$, $y$ and $z$ are odds, the numbers $X$ are also odd and, therefore, they can be represented as differences of two squares $([X+1]/2)^2$ and $([X-1]/2)^2$. Taking such an $X$ and $k \in \mathbb{Z}^+, \mu \in \mathbb{Z}^+, \nu \in \mathbb{Z}^+$ so that $X+1 = 2kz^{\mu(r-2)}$ and $X-1 = 2\nu$, we see that $y^{\sigma q-2} \left\{ [kz^{\mu(r-2)}]^2 - \nu^2 \right\} = \pm[z^2 \mp \tau_0 y^2] \mp \pi$ for a prime $\pi$ and, therefore, if $m^* = kz^{\mu(r-2)}$, the integer $n^* = -kz^{\mu(r-2)} \pm \nu$ will be a root of the equation $(y^{\sigma q-2})n^{*2} + (2m^* y^{\sigma q-2})n^* + [\pm(z^2 \mp \tau_0 y^2) \mp \pi] = 0$. We infer $(\pm y^{\sigma q-2})n^{*2} + (\pm 2m^* y^{\sigma q-2})n^* + (z^2 \mp \tau_0 y^2) - \pi = z^2 \mp y^2 \tau_{m^*,n^*} - \pi = 0$, and, thus, the proof of Proposition 1 is completed.

**Proposition 2**. *We have $\tau_{m^*,n^*} z^{r-2} - y^{q-2} > 1$ and $gcd(x^{p-2}, \tau_{m^*,n^*} z^{r-2} - y^{q-2}) = 1$.*

*Proof.* Let $m^*$ and $n^*$ be as in Proposition 1. Because of this choice, it is easily verified that $\tau_{m^*,n^*} z^{r-2} - y^{q-2} > 1$. For the rest of the proof. To get a contradiction, suppose that $gcd(x^{p-2}, \tau_{m^*,n^*} z^{r-2} - y^{q-2}) > 1$. This means that there were a natural number $\varrho > 1$ such that $x^{p-2} = \varrho c_1$ and $\tau_{m^*,n^*} z^{r-2} - y^{q-2} = \varrho c_3$ for some integers $c_1$ and $c_3$. Multiplication of the first equation by $x^2$, gives

$$x^p = \varrho(c_1 x^2), \qquad (1)$$

while multiplication of the second one by $y^2$ and $z^2$ gives

$$y^2 \tau_{m^*,n^*} z^{r-2} - y^q = \varrho(c_3 y^2) \Rightarrow \pm y^q = \pm y^2 \tau_{m^*,n^*} z^{r-2} \mp \varrho(c_3 y^2) \qquad (2)$$

and

$$z^2 \tau_{m^*,n^*} z^{r-2} - z^2 y^{q-2} = \varrho(c_3 z^2) \Rightarrow \tau_{m^*,n^*} z^r = z^2 y^{q-2} + \varrho(c_3 z^2), \qquad (3)$$

respectively. Adding (2) to (1), we obtain $x^p \pm y^q = \varrho(c_1 x^2 \mp c_3 y^2) \pm y^2 \tau_{m^*,n^*} z^{r-2}$ or, by hypothesis, $z^r = \varrho(c_1 x^2 \mp c_3 y^2) \pm y^2 \tau_{m^*,n^*} z^{r-2}$ which can equivalently be written as

$$z^{r-2}(z^2 \mp y^2 \tau_{m^*,n^*}) = \varrho(c_1 x^2 \mp c_3 y^2). \qquad (4)$$

Similarly, subtracting (3) from the multiple of (1) by $\tau_{m^*,n^*}$, we get $\tau_{m^*,n^*}(x^p - z^r) = \varrho(\tau_{m^*,n^*} c_1 x^2 - c_3 z^2) - z^2 y^{q-2}$ or, by hypothesis, $\mp \tau_{m^*,n^*} y^q = \varrho(\tau_{m^*,n^*} c_1 x^2 - c_3 z^2) - z^2 y^{q-2}$ which can equivalently be written as

$$y^{q-2}(z^2 \mp \tau_{m^*,n^*} y^2) = \varrho(\tau_{m^*,n^*} c_1 x^2 - c_3 z^2). \qquad (5)$$

Having regard to the Proposition 1, from the relations (4) and (5), it follows that $\varrho > 1$ divides both $z^{r-2}$ and $y^{q-2}$, which contradicts the statement that $gcd(x, y, z) = 1$ under the condition that $x^p \pm y^q = z^r$. So, it is impossible to hold $gcd(x^{p-2}, \tau_{m^*,n^*} z^{r-2} - y^{q-2}) > 1$. Hence, $gcd(x^{p-2}, \tau_{m^*,n^*} z^{r-2} - y^{q-2}) = 1$, and the proof of Proposition 2 is complete.

We are now in position to give the main result of the paper.

**Theorem.** *There is no non-trivial solution to the generalized Fermat equation $x^p + y^q = z^r$, once we assume that $\min\{p, q, r\} \geq 3$ and $\gcd(x, y, z) = 1$.*

*Proof.* Let $m^*$ and $n^*$ be as in Proposition 1. Let also $U, V, \tilde{U}, \tilde{V} \in \mathbb{Z} \setminus \{0\}$ with $U\tilde{V} - \tilde{U}V \neq 0$. Then, the system

$$\begin{bmatrix} (U)X + (-V)Y + (\tau_{m^*,n^*} V)Z = \pm 1 \\ (\tilde{U})X + (-\tilde{V})Y + (\tau_{m^*,n^*} \tilde{V})Z = \pm 1 \\ (x^2)X + (\pm y^2)Y + (-z^2)Z = 0 \end{bmatrix}$$

with determinant $(U\tilde{V} - \tilde{U}V)(z^2 \mp y^2\tau_{m^*,n^*}) \neq 0$ has unique solution

$$X = \mp \frac{V - \tilde{V}}{U\tilde{V} - \tilde{U}V}, Y = \mp \frac{(U-\tilde{U})z^2 + (V-\tilde{V})\tau_{m^*,n^*} x^2}{(U\tilde{V} - \tilde{U}V)(z^2 \mp y^2\tau_{m^*,n^*})} \text{ and } Z = -\frac{(U-\tilde{U})y^2 \pm (V-\tilde{V})x^2}{(U\tilde{V} - \tilde{U}V)(z^2 \mp y^2\tau_{m^*,n^*})}.$$

From Proposition 2, it holds true $\gcd(x^{p-2}, \tau_{m^*,n^*} z^{r-2} - y^{q-2}) = 1$, so an application of Bezout's Identity shows that the equation $(x^{p-2})u + (\tau_{m^*,n^*} z^{r-2})v + (-y^{q-2})v = \pm 1$ can be solved in $\mathbb{Z}^2$. Given any partial integer solution $(u, v) \in \mathbb{Z}^2$ of this equation and $\varkappa, \ell \in \mathbb{Z}, \varkappa \neq \ell$, let us take the induced solutions $U, V, \tilde{U}, \tilde{V} \in \mathbb{Z} \setminus \{0\}$ defined by

$$U := u + \varkappa(\tau_{m^*,n^*} z^{r-2} - y^{q-2}), V := v - \varkappa x^{p-2},$$
$$\tilde{U} := u + \ell(\tau_{m^*,n^*} z^{r-2} - y^{q-2}) \text{ and } \tilde{V} := v - \ell x^{p-2}.$$

Observe that $U, V, \tilde{U}, \tilde{V} \neq 0$. Further, since $(U\tilde{V} - \tilde{U}V) = \pm(\varkappa - \ell) \neq 0$ and it is easily seen that the unique solution of the system is exactly $(X, Y, Z) = (x^{p-2}, y^{q-2}, z^{r-2})$. But, on the other hand, considering the above-mentioned expression of the unique solution, we get $X = x^{p-2}$, $Y = (-z^2)D + (x^p\tau_{m^*,n^*})E$ and $Z = (\mp y^2)D + (x^p)Ex^p$, where we have used the notation

$$D := \frac{\tau_{m^*,n^*} z^{r-2} - y^{q-2}}{z^2 \mp y^2\tau_{m^*,n^*}} \text{ and } E = \frac{1}{z^2 \mp y^2\tau_{m^*,n^*}}.$$

Thus, especially in this case, with this formulation, we found

$$\begin{bmatrix} y^{q-2} = (-z^2)D + (x^p\tau_{m^*,n^*})E \\ z^{r-2} = (\mp y^2)D + (x^p)E \end{bmatrix}.$$

Since $x^p \pm y^q = z^r$, this formulation guarantees that the system of equations

$$\begin{bmatrix} (1 - z^2E)Z \pm (y^2E)Y = \mp y^2D \\ (\pm z^2\tau_{m^*,n^*}E)Z + (-y^2\tau_{m^*,n^*}E \mp 1)Y = \pm z^2D \end{bmatrix}$$

would have a solution at $(Z, Y) = (z^{r-2}, y^{q-2})$, which is impossible, since the determinant of the system equals zero and both lines representing it must in any case not be identical, otherwise we would have $\tau_{m^*,n^*} = 0$ and $z^2 = -y^2$. We therefore conclude that there is no non-trivial solution to the generalized Fermat equation $x^p \pm y^p = z^p$ and, thus, the proof is complete.

We finish by giving two immediate consequences.

**Corollary 1** (Fermat's Last Theorem) *There is no non-trivial integer solution to the Fermat Equation $x^p + y^p = z^p$ provided that $p \geq 3$ and $\gcd(x, y, z) = 1$.*

**Corollary 2** (Irrationality of Binomial Expansions) *If $p, q, r \geq 3$, there is no $(x, y) \in \mathbb{Z}^2$ with $\gcd(x, y) = 1$ and satisfying $(x^p + y^q)^{1/r} \in \mathbb{Q}$.*

*Proof.* Suppose there is a rational number $(P/Q) \in \mathbb{Q}$, with $P, Q \in \mathbb{Z}$ and $gcd(P, Q) = gcd(x, y, P) = 1$, satisfying $x^p + y^q = (P/Q)^r$. This implies that $P^r = Q^r(x^p + y^q)$. Since $gcd(P, Q) = 1$, we infer $P^r / (x^p + y^q)$, which means that there is a $k \in \mathbb{Z}$ satisfying $(x^p + y^q) = kP^r$. It follows that $P^r = Q^r k P^r$, which guarantees that $Q^r k = 1$. But, this only applies if $Q = k = 1$. In such a case, we would have $x^p + y^q = P^r$ and also $gcd(x, y, P) = 1$, which, in view of the Theorem above, is impossible to hold true.

## References

[1]   M. A. Bennett, I. Chen, S. R. Dahmen and S. Yazdani: *Generalized Fermat equations: a miscellany*, Int. J. Number Theory **11** (2015), 1–28

[2]   H. Darmon and A. Granville: *On the equations $z^m = F(x, y)$ and $Ax^p + By^q = Cz^r$*, Bull.London Math. Soc. Textbf **27** (1995), 513–543

[3]   Michiel Hazewinkel, ed. (2001) [1994]: *Hilbert problems*, **Encyclopedia of Mathematics**, Springer Science + Business Media B.V. / Kluwer Academic Publishers, ISBN 978-1-55608-010-4. For an electronic version, see https://www.encyclopediaofmath.org/index.php/Hilbert_problems

[4]   David Hilbert: *Mathematical Problems*, Bulletin of the American Mathematical Society **8(10)** (1902), pp. 437-479. (Also, Bulletin (New Series) of the American Mathematical Society **37(4)**, pp. 407-436, S 0273-0979(00)00881-8. Article electronically published on June 26, 2000. Text on the web: http://www.ams.org/journals/bull/2000-37-04/S0273-0979-00-00881-8/S0273-0979-00-00881-8.pdf ) Earlier publications (in the original German) appeared in Göttinger Nachrichten, 1900, pp. 253-297, and Archiv der Mathematik und Physik, 3dser., vol. 1 (1901), pp. 44-63, 213-237

[5]   R. D. Mauldin: *A Generalization of Fermat's Last Theorem: The Beal Conjecture and Prize Problem*, Notices of the AMS **44** (11) (1997), 1436–1439

[6]   P. Mihăilescu: *Primary Cyclotomic Units and a Proof of Catalan's Conjecture*, J. Reine Angew. Math.**572** (2004), 167–195

[7]   J.F. Nash and M.Th. Rassias(eds.): Open Problems in Mathematics, Springer (2016), 543 pages

[8]   A. Wiles: *Modular elliptic curves and Fermat's Last Theorem*, Ann. Math. **141** (1995),443–551